\documentstyle[12pt]{article}
\begin{document}
\parindent 16pt
\title{\Large\bf Non-Conflicting Ordering Cones and Vector Optimization in Inductive Limits\thanks{This research was supported by the National Natural
Science Foundation of China (10871141).} }
\author { Jing-Hui Qiu\\
 {\small \sl School of Mathematical Sciences, Soochow University,  Suzhou 215006,  China}\\
{\small\sl  E-mail: qjhsd@sina.com, jhqiu@suda.edu.cn }}
\date{}
\maketitle
\begin{center}
\begin{minipage}{124mm}
\vskip 0.3cm {{\bf Abstract}\ \  Let $(E,\xi)={\rm ind}(E_n, \xi_n)$
be an inductive limit of a sequence $(E_n, \xi_n)_{n\in N}$ of
locally convex spaces and let every step $(E_n, \xi_n)$ be endowed
with a partial order by a pointed convex (solid)  cone $S_n$. In the
framework of inductive limits of partially ordered locally convex
spaces, the notions of lastingly efficient points, lastingly weakly
efficient points and lastingly globally properly efficient points
are introduced. For several ordering cones, the notion of
non-conflict is introduced. Under the requirement that the sequence
$(S_n)_{n\in N}$ of ordering cones is non-conflicting, an existence
theorem on lastingly weakly efficient points is presented. From
this, an existence theorem on lastingly globally properly efficient
points is deduced. \\

}\\

{\bf MR(2000) Subject Classification} 46A03, 46A13, 90C48\\

\end{minipage}
\end{center}
\vskip 1cm \baselineskip 20pt

\section* {\large\bf 1 \ \  Introduction}

\ \ \ Let $(X, \tau)$ be a real locally convex Hausdorff
topological vector space (briefly, denoted by a locally convex
space) and $X^*$ be its topological dual. A set $S\subset X$ is
said to be a cone if  $\alpha s\in S$ for any $s\in S$ and any
$\alpha\geq 0$; and a convex cone if in addition $S+S\subset S$. A
cone $S$ is said to be pointed if $S\cap (-S)=\{0\}$. Furthermore,
${\rm int}_{\tau} S$ (briefly, denoted by ${\rm int}\, S$ if no
confusions) denotes the topological interior of $S$ in $(X, \tau
)$ and $S^+$ denotes the dual cone of $S$, i.e.,
$$S^+\,=\,\{f\in X^*:\, f(s)\geq 0,\ \forall s\in S\}.$$ If
${\rm int}_{\tau} S\not=\emptyset$, we call $S$ a solid cone in
$(X,\tau)$. As is well-known (for example, see [1, 2, 3]), a pointed
convex cone $S$ specifies a partial order in $X$ as follows:
$$ {\rm for}\ \  x, y\in X,\ \ x\leq_S y\ \ \ {\rm iff}\ \ \ y-x\in S.$$ Obviously,  the partial order $\leq_S$ determined
by a pointed convex cone $S$ satisfies the following properties:

(i) $x\leq_S x$;

(ii) if $x\leq_S y$ and $y\leq_S z$, then $x\leq_S z$;

(iii) if $x\leq_S y$ and $y\leq_S x$, then $x=y$.

In this case, the pointed convex cone $S$ is called an ordering
cone. Let $A\subset X$ be nonempty and $S\subset X$ be an ordering
cone. We denote $E(A, S)$ the set of efficient points (i.e., Pareto
minimal points) of $A$ with respect to the ordering cone $S$; that
is,
$$a_0\in E(A, S)\ \ \  {\rm iff}\ \ \ (A-a_0)\cap
(-S)\,=\,\{0\}.$$ Moreover, if ${\rm int}\,S\not=\emptyset$, i.e.,
$S$ is a pointed convex solid cone, we denote $w-E(A,S)$ the set of
weakly efficient points (i.e., weakly Pareto minimal points) of $A$
with respect to the ordering cone $S$; that is,
$$a_0\in w-E(A, S)\ \ \ {\rm iff}\ \ \ (A-a_0)\cap
(-{\rm int}\,S)\,=\,\emptyset.$$ As observed by Kuhn and Tucker [4]
and later by Geoffrion [5], some efficient points exhibit certain
abnormal properties. To eliminate such anomalous efficient points,
various concepts of proper efficiency have been introduced (for
example, see [6, 7, 8]). Efficiency, weak efficiency and proper
efficiency in locally convex spaces have already been investigated;
for example, see [1-8] and the references therein. In this paper, we
consider vector optimization in the framework of inductive limits of
partially ordered locally convex spaces. We shall introduce new
concepts of efficiency, weak efficiency and proper efficiency in the
new framework and give some existence results. First we recall some
concepts concerning inductive limits (for example, see [9-11). Let
$(E_1,\xi_1)\subset (E_2, \xi_2)\subset\cdots $ be an increasing
sequence of locally convex spaces with continuous inclusions $(E_n,
\xi_n)\,\rightarrow\, (E_{n+1}, \xi_{n+1})$ for every $n$. Here, we
require that each $E_n$ is properly included in $E_{n+1}$. If the
union $E=\cup_{n=1}^{\infty} E_n$ is endowed with the finest locally
convex topology $\xi$ such that all inclusions $i_n:\, (E_n,
\xi_n)\,\rightarrow\, E$ are continuous, then $(E,\xi)$ is called
the locally convex inductive limit of a sequence $(E_n, \xi_n)_{n\in
N}$ of locally convex spaces; and it is denoted by ${\rm ind}(E_n,
\xi_n)$. Also, every $(E_n, \xi_n)$ is called a step of the
inductive limit ${\rm ind}(E_n, \xi_n)$. In general, $(E,\xi)$
needn't be Hausdorff even though every $(E_n, \xi_n)$ is Hausdorff
(see [10]). But in this paper we always assume that $(E, \xi)$ is
Hausdorff. For each $n\in N$, let $S_n$ be a pointed convex cone in
$(E_n, \xi_n)$, which specifies a partial order $\leq_{S_n}$ in
$E_n$. Thus, we have a sequence of partially ordered locally convex
spaces, i.e., $(E_n, \xi_n, \leq_{S_n})_{n\in N}$. Suppose that one
is placed in the first step $(E_1, \xi_1)$ and $A\subset E_1$ is a
feasible set. Certainly, one is concerned about the set of efficient
points  of $A$ with respect to the ordering cone $S_1$, i.e., $E(A,
S_1)$. An interesting problem is: Does there exist a point $a_0\in
A$ such that not only $a_0\in E(A, S_1)$ but also $a_0\in E(A, S_n)$
for all $n\in N$?

Moreover, assume that every ${\rm int}_{\xi_n}
S_n\,\not=\,\emptyset$, where ${\rm int}_{\xi_n} S_n$ denotes the
topological interior of $S_n$ in $(E_n, \xi_n)$. If $a_0\in A$ such
that $(A-a_0)\cap(-{\rm int}_{\xi_n} S_n) =\emptyset$, we denote
$a_0\in w-E_n(A, S_n)$. Here, ${\rm int}_{\xi_n} S_n$ is related to
the topology $\xi_n$ of $E_n$, so we denote the weakly efficient
point set  by $w-E_n(A, S_n)$ rather than $w-E(A, S_n)$. Similarly,
an interesting problem is: Does there exist a point $a_0\in A$ such
that not only $a_0\in w-E_1(A, S_1)$, but also $a_0\in w-E_n(A,
S_n)$ for all $n\in N$?

For convenience, we introduce the following notions.\\

{\bf Definition 1.1} \ \ Let $(E,\xi)\,=\,{\rm ind}(E_n, \xi_n)$ be
an inductive limit of a sequence $(E_n, \xi_n)_{n\in N}$ of locally
convex spaces and $S_n$ be an ordering cone in $E_n$ for every $n$.
Suppose that $A\subset E_1$ is nonempty and $a_0\in A$. Then $a_0$
is said to be a lastingly efficient point of $A$ with respect to the
sequence $(S_n)_{n\in N}$ of ordering cones if
$$a_0\,\in\,\bigcap_{n=1}^{\infty} E(A, S_n).$$
Moreover, suppose that ${\rm int}_{\xi_n} S_n\,\not=\,\emptyset$
for every $n$. Then $a_0\in A$ is said to be a lastingly weakly
efficient point of $A$ with respect to the sequence $(S_n)_{n\in
N}$ of ordering cones if
$$a_0\,\in\,\bigcap_{n=1}^{\infty} w-E_n(A, S_n).$$

In terms of Definition 1.1, the above problems become: Does there
exist a point $a_0\in A$ such that $a_0$ is a lastingly efficient
point (resp. lastingly weakly efficient point) with respect to the
sequence $(S_n)_{n\in N}$ of ordering cones?

In next section, we shall give an existence theorem on lastingly
weakly efficient points.\\

\section* {\large\bf 2 \ \  Non-conflicting ordering cones and existence of lastingly weakly efficient points}

\ \ \ \ \ In order to obtain a satisfying existence result of
lastingly weakly efficient points, the following ``non-conflict"
requirement on the sequence $(S_n)_{n\in N}$ of ordering cones seems
to be necessary by one's experience.\\

{\bf Definition 2.1} \ \   Let $(X,\tau )$ be a locally convex space
and $S_1,\ S_2,\ \cdots,\ S_n$ be ordering cones in $(X,\tau )$. The
ordering cones $S_1,\ S_2,\ \cdots,\ S_n$ are called non-conflicting
if
$$-S_j\,\not\subset\,{\rm cl}_{\tau}(S_1+S_2+\cdot +S_n),\ \ \ {\rm
for}\  j=1,2, \cdots, n.$$ Here ${\rm cl}_{\tau}(S_1+S_2+\cdots
+S_n)$ denotes the closure of $S_1+S_2+\cdots+S_n$ in $(X,\tau )$.\\

{\bf Definition 2.2} \ \  let $(E,\xi)\,=\,{\rm ind}(E_n,\xi_n)$ and
$S_n$ be a pointed convex cone in $(E_n,\xi_n)$ for every $n\in N$.
The sequence $(S_n)_{n\in N}$ is said to be non-conflicting in ${\rm
ind}(E_n,\xi_n)$ if for each $n\in N$, $S_1,\ S_2,\ \cdots,\ S_n$
are non-conflicting in $(E_n,\xi_n)$. That is, for each $n\in N$,
$$-S_j\,\not\subset\,{\rm cl}_{\xi_n}(S_1+S_2+\cdots
+S_n) \ \ \ {\rm for}\ j=1,2,\cdots,n.$$ Here  ${\rm
cl}_{\xi_n}(S_1+S_2+\cdots +S_n)$ denotes the closure of
$S_1+S_2+\cdots +S_n$ in $(E_n, \xi_n)$.\\

{\bf Lemma 2.1} \ \ {\sl Let $(X, \tau )$ be a locally convex space
and $S\subset X$ be a pointed convex cone with ${\rm int}\,
S\,\not=\,\emptyset$. If $A\subset X$ is weakly countably compact in
$(X, \tau)$, then:

{\rm (i)} there exists $f\in S^+\backslash\{0\}$ and $a_0\in A$ such
that
$$f(a_0)\,=\,{\rm inf}\{f(a):\, a\in A\};$$

{\rm (ii)}  $w-E(A, S)\,\not=\,\emptyset$.}\\

{\bf Proof}\ \  Since $S$ is a pointed convex cone, $0\not\in {\rm
int}\,S$. Hence there exists $f\in X^*\backslash\{0\}$ such that
$$f({\rm int}\,S)>0. \eqno{(1)}$$
Here and in the following, for any nonempty subset $B$ of $X$,
$f(B)\,>\,0$ means that $f(x)>0,\  \forall x\in B$. Similarly,
$f(B)\,\geq\,0$ means that $f(x)\,\geq\,0,\  \forall x\in B$. Since
$A$ is weakly countably compact in $(X, \tau)$, we know that $A$ is
bounded in $(X, \tau)$ and
$$-\infty\,<\,\alpha={\rm inf}\{f(a):\,a\in A\}\,<\,+\infty.$$
Take a sequence $(a_n)\subset A$ such that
$$f(a_n)\,\rightarrow\,\alpha,\ \ \ n\rightarrow \infty.
\eqno{(2)}$$ Since $A$ is weakly countably compact, $(a_n)$ has a
weakly cluster point $a_0\in A$. Thus $f(a_0)$ is a cluster point of
$(f(a_n))$ in the real line $R$. Hence there exists a sequence
$n_1<n_2<\cdots$ of natural numbers such that
$$f(a_{n_i})\,\rightarrow\,f(a_0),\ \ \  i\rightarrow\infty. \eqno{(3)}$$
By (2) and (3), we have $f(a_0)=\alpha$. Thus (i) is shown. As
$f(a_0)=\inf \{f(a):\,a\in A\}$, we have
$$f(A-a_0)\,\geq\,0. \eqno{(4)}$$
By (1) and (4), we have
$$(A-a_0)\cap (-{\rm int}\, S)\,=\,\emptyset.$$
Therefore $a_0\in w-E(A, S)$ and (ii) has been proven.
\hfill\framebox[2mm]{}\\

{\bf Lemma 2.2} \ \ {\sl Let $(X, \tau)$ be a locally convex space,
$S\subset X$ be a pointed convex cone with ${\rm
int}\,S\not=\emptyset$ and $A\subset X$ be nonempty. If $w-E(A,
S)\not=\emptyset$, then $w-E(A, S)$ is a closed subset of
$A$.}\\

{\bf Proof}\ \  Assume that a net $(x_{\delta})_{\delta\in \Delta}$
in $w-E(A, S)$ is convergent to $x_0\in A$. Clearly, for each
$\delta\in\Delta$,
$$(A-x_{\delta})\cap (-{\rm int}\,S)\,=\,\emptyset.\eqno{(5)}$$
We are going to prove that
$$(A-x_0)\cap (-{\rm int}\, S)\,=\,\emptyset.\eqno{(6)}$$
If not, there exists $a\in A$ such that
$$a-x_0\in -{\rm int}\,S.$$
Since $a-x_{\delta}\,\rightarrow\,a-x_0$ and $-{\rm int}\,S$ is an
open neighborhood of $a-x_0$, there exists $\delta_0\in\Delta$ such
that
$$a-x_{\delta}\in -{\rm int}\,S\ \ \ {\rm for\ \ all}\ \
\delta\geq\delta_0.$$ This contradicts (5). Thus we have shown (6)
and $x_0\in w-E(A, S)$. This means that $w-E(A,S)$ is a closed
subset of $A$.\hfill\framebox[2mm]{}\\

Now we present our main result as follows.\\

{\bf Theorem 2.1} \ \ {\sl Let $(E,\xi)\,=\,{\rm ind}(E_n, \xi_n)$
be an inductive limit of a sequence $(E_n, \xi_n)_{n\in N}$ of
locally convex spaces. For each $n\in N$, let $S_n$ be a pointed
convex cone in $(E_n, \xi_n)$ with ${\rm int}_{\xi_n}
S_n\,\not=\,\emptyset$. Moreover, assume that the sequence
$(S_n)_{n\in N}$ is non-conflicting. If $A\subset E_1$ is countably
compact in $(E_1, \xi_1)$, then the set of lastingly weakly
efficient points of $A$ with respect to the sequence $(S_n)_{n\in
N}$ is nonempty, i.e.,
$$\bigcap\limits_{n=1}^{\infty} w-E(A,
S_n)\,\not=\,\emptyset.$$}\\

{\bf Proof}\ \  By the assumption, $(S_n)_{n\in N}$ is
non-conflicting. Hence for every $n\in N$, we have
$$-S_j\,\not\subset\,{\rm cl}_{\xi_n}(S_1+S_2+\cdots+S_n),\ \ \
{\rm for}\ \ j=1,2,\cdots, n.$$ Thus, for each $j\in\{1,2,\cdots,
n\}$, there exists $y_j\in S_j$ such that $-y_j\,\not\in\,{\rm
cl}_{\xi_n}(S_1+S_2+\cdots+ S_n)$. By the Hahn-Banach separation
theorem, there exists $f_j\,\in\,(E_n,\xi_n)^*\backslash\{0\}$ and
$\eta_j > 0$ such that
$$f_j(-y_j)\,\leq\,-\eta_j\,<\,0\,\leq\,f_j(S_1+S_2 +\cdots
S_n).\eqno{(7)}$$ Obviously,
$$f_j(S_i)\,\geq\,0,\ \ \ {\rm for}\ i=1,2,\cdots,n;\
j=1,2,\cdots,n.$$ Put $f:=\,\sum_{j=1}^{n}f_j$. Then $f\in (E_n,
\xi_n)^*\backslash\{0\}$ and $f(S_i)\,\geq\,0$, for $i=1,2,\cdots,
n$. Next we show that
$$f({\rm int}_{\xi_i}S_i)\,>\,0,\ \ \ {\rm for}\ i=1, 2,\cdots,
n.\eqno{(8)}$$ If not, there exists $i$ with $1\leq i\leq n$ and
there exists $y^{\prime}\,\in\,{\rm int}_{\xi_i}S_i$ such that
$f(y^{\prime})\,=\,\sum_{j=1}^n f_j(y^{\prime})\,=\,0$. Since every
$f_j(y^{\prime})\,\geq\,0$, we have $f_j(y^{\prime})\,=\,0$, for
$j=1,2,\cdots, n.$ Particularly, $f_i(y^{\prime})\,=\,0$. Since
$y^{\prime}\,\in\,{\rm int}_{\xi_i}S_i$, there exists an absolutely
convex 0-neighborhood $U_i$ in $(E_i, \xi_i)$ such that $y^{\prime}
+ U_i\,\subset\, S_i$. Remark that $f_i(S_i)\,\geq\,0$, so we have
$f_i(U_i)\,=\,f_i(y^{\prime} +U_i)\,\geq\,0$. Since $U_i$ is an
absolutely convex set, we also have that $f_i(U_i)\,\leq\,0$. Thus,
$$f_i(U_i)\,=\,\{0\}\ \ {\rm and\ \ hence}\ \
f_i(E_i)\,=\,\{0\}.\eqno{(9)}$$ On the other hand, from (7) we
have
$$f_i(y_i)\,\geq\,\eta_i\,>\,0\ \ \ {\rm and}\ \ \ y_i\in
S_i\subset E_i,$$ which contradicts (9). Thus we have shown that (8)
is true. Since $A$ is countably compact in $(E_1,\xi_1)$ and the
inclusion $(E_1, \xi_1)\,\rightarrow\,(E_n, \xi_n)$ is continuous,
we know that $A$ is also countably compact in $(E_n, \xi_n)$ and
hence weakly countably compact in $(E_n, \xi_n)$. Thus, by Lemma 2.1
there exists a point $a_0\in A$ such that
$$f(a_0)\,=\,{\rm inf}\{f(a):\,a\in A\},\ \ {\rm i.e.,}\ \
f(A-a_0)\geq 0.\eqno{(10)}$$ Combining (8) and (10), we have
$$(A-a_0)\,\cap\,(-{\rm int}_{\xi_i} S_i)\,=\,\emptyset,\ \ \
i=1,2,\cdots,n.$$ That is,
$$a_0\in \bigcap_{i=1}^n w-E_i(A, S_i).$$
By Lemma 2.2, $w-E_i(A, S_i)$ is a $\xi_i$-closed subset of $A$ and
it is also a $\xi_1$-closed subset of $A$ since the topology $\xi_1$
is finer than one induced by $\xi_i$. Thus, for every $n\in N$,
$$\bigcap_{i=1}^n w-E_i(A, S_i)$$
is a nonempty $\xi_1$-closed subset of $A$.
 Since $A$ is countably compact in $(E_1, \xi_1)$, we conclude
 that
 $$\bigcap_{i=1}^{\infty} w-E_i(A, S_i)\,\not=\,\emptyset.$$
 Thus the proof is completed. \hfill\framebox[2mm]{}\\

\section* {\large\bf 3 \ \  Existence of lastingly globally properly efficient points}

\ \ \ \ As we know, there are various concepts of proper efficiency
in vector optimization, for example, see [6, 7, 8]. In the framework
of inductive limits, we may consider various concepts of lastingly
proper efficiency. In this section, we introduce the notion of
lastingly globally properly efficient points (in fact, i.e.,
lastingly generalized Henig properly efficient points) and
investigate the existence of this kind of properly efficient points.
First we recall some related concepts.\\

{\bf Definition 3.1}(see [7])\ \  Let $X$ be a locally convex space,
$A\subset X$ be a nonempty set and $S\subset X$ be an ordering cone.
A point $a\in A$ is said to be a globally properly efficient point
of $A$, denoted by $a\in GHe(A, S)$, iff there exists a convex cone
$W$ such that $S\backslash\{0\}\subset {\rm int}W$ and $a\in E(A,
W)$. Here, $W$ is said to be a dilating cone of $S$.\\

{\bf Definition 3.2}(see [7])\ \  Let $X$, $A$ and $S$ be the same
as in Definition 3.1. Moreover, let $S$ have a base $\Theta$ (i.e.,
$\Theta\subset S$ is a convex set such that $S={\rm cone}(\Theta)$
and $0\not\in {\rm cl}(\Theta)$). Clearly, there exists a convex
0-neighborhood $U$ such that $0\not\in \Theta +U$. A point $a\in A$
is said to be a Henig properly efficient point of $A$ with respect
to $\Theta$, denoted by $a\in HE(A, \Theta)$, iff there exists a
convex 0-neighborhood $V\subset U$ such that
$${\rm cl}{\rm cone}(A-a)\cap -{\rm cone}(\Theta +V)\,=\,\{0\}.$$
Let $B(S)$ denote the family of all bases of $S$. A point $a\in A$
is said to be a Henig properly efficient point (respectively, a
generalized Henig properly efficient point) with respect to $S$,
denoted by $a\in HE(A, S)$ (respectively, $a\in GHE(A, S)$), if
$a\in \bigcap_{\Theta\in B(S)}HE(a, \Theta)$ (respectively, $a\in
\bigcup_{\Theta\in B(S)} HE(A, \Theta)$).\\

Liu and Song (see [8, Theorem 3.4])  showed that $GHE(A, S) \subset
GHe(A,S)$. Moreover, the author and Hao (see [12, Theorem 4.2])
proved the following.\\

{\bf Lemma 3.1}\ \  {\sl $GHE(A, S) =GHe(A, S)$.}\\

  In the
framework of inductive limits of partially ordered locally convex
spaces, we give the definition of lastingly globally properly
efficient points (i.e., lastingly generalized Henig properly
efficient points) as
follows.\\

{\bf Definition 3.3}\ \  Let $(E, \xi) ={\rm ind}(E_n, \xi_n)$ be an
inductive limit and $S_n$ be an ordering cone in $E_n$ for every
$n\in N$. Suppose that $A\subset E_1$ is nonempty and $a_0\in A$.
Then $a_0$ is said to be a globally properly efficient point in the
partially ordered locally convex space $(E_n, \xi_n, S_n)$, denoted
by $a_0\in GHe_n(A, S_n)$, iff there exists a convex cone $W_n$ in
$E_n$ such that $S_n\backslash\{0\}\subset {\rm int}_{\xi_n} W_n$
and $a_0\in E(A, W_n)$. Moreover, $a_0$ is said to be a lastingly
globally properly efficient point of $A$ with respect to the
sequence $(S_n)_{n\in N}$ of ordering cones iff $a_0\in
\bigcap_{n=1}^{\infty} GHe_n(A, S_n)$.\\

{\bf Definition 3.4}\ \  Let $(E,\xi)={\rm ind}(E_n, \xi_n)$, $S_n$,
$A$ and $a_0$ be the same as in Definition 3.3. Then $a_0$ is said
to be a generalized Henig properly efficient point in the partially
ordered locally convex space $(E_n, \xi_n, S_n)$, denoted by $a_0\in
GHE_n(A, S_n)$, iff there exists a base $\Theta_n$ of $S_n$ and a
convex 0-neighborhood $V_n$ in $(E_n, \xi_n)$ such that
$${\rm cl} {\rm cone}(A-a_0)\cap -{\rm
cone}(\Theta_n+V_n)\,=\,\{0\}.$$ Moreover, $a_0$ is said to be a
lastingly generalized Henig properly efficient point of $A$ with
respect to the sequence $(S_n)_{n\in N}$ of ordering cones iff
$a_0\in \bigcap_{n=1}^{\infty} GHE_n(A, S_n)$.\\

By Lemma 3.1, we see that the set of lastingly globally properly
efficient points and the set of lastingly generalized Henig properly
efficient points are the same. From Theorem 2.1, we deduce an
existence theorem on lastingly globally properly efficient points
(i.e., lastingly generalized Henig properly efficient points) as
follows.\\

{\bf Theorem 3.1}\ \ {\sl Let $(E, \xi)={\rm ind}(E_n, \xi_n)$ be an
inductive limit of a sequence $(E_n, \xi_n)_{n\in N}$ of locally
convex spaces. For each $n\in N$, let $S_n$ be a pointed convex cone
and $W_n$ be a dilating cone of $S_n$ in $(E_n, \xi_n)$. Moreover,
assume that the sequence $(W_n)_{n\in N}$ is non-conflicting. If
$A\subset E_1$ is countably compact in $(E_1, \xi_1)$, then the set
of lastingly globally properly efficient points of $A$ with respect
to the sequence $(S_n)_{n\in N}$ is nonempty, i.e.,
$\bigcap_{n=1}^{\infty} GHe_n(A, S_n)\not=\emptyset$. Equivalently,
we have $\bigcap_{n=1}^{\infty} GHE_n(A, S_n)\not=\emptyset$.}\\

{\bf Proof}\ \  Without loss of generality, we may assume that every
$W_n$ is a pointed cone. Or else, we may take $W_n^{\prime}$, where
$W_n^{\prime}={\rm int}_{\xi_n}W_n \cup\{0\}$. It is easy to show
that $(W^{\prime}_n)_{n\in N}$ is non-conflicting provided that
$(W_n)_{n\in N}$ is non-conflicting. Thus, applying Theorem 2.1, we
conclude that $\bigcap_{n=1}^{\infty} w-E_n(A, W_n)\not=\emptyset$.
Let $a_0\in \bigcap_{n=1}^{\infty}w-E_n(A, W_n)$. Then $a_0\in
\bigcap_{n=1}^{\infty} GHe_n(A, S_n)$. Equivalently,
$a_0\in\bigcap_{n=1}^{\infty} GHE_n(A, S_n)$. Hence
$$\bigcap_{n=1}^{\infty} GHe_n(A, S_n) = \bigcap_{n=1}^{\infty}
GHE_n (A, S_n)\not=\emptyset.$$\hfill\framebox[2mm]{}

\noindent {\bf References} \vskip 12pt
\begin{description}
\item{[1]} D. T. Luc, Theory of Vector Optimization, Berlin,
 Springer, Berlin, 1989.
\item{[2]} J. Jahn, Vector Optimization, Springer, Berlin, 2004.
\item{[3]} G. Y. Chen, X. X. Huang, X. Q. Yang, Vector Optimization
-- Set-Valued and Variational Analysis, Springer, Berlin, 2005.
\item{[4]} H. Kuhn, A. Tucker, Nonlinear programming, in:
Proceedings of the Second Berkeley Symposium on Mathematical
Statistics and Probability, University of California Press,
Berkeley, CA, 1951, pp. 481-492.
\item{[5]} A. M. Geoffrion, Proper efficiency and the theory of
vector maximization, J. Math. Anal. Appl. {\bf 22} (1968) 618-630.
\item{[6]} A. Guerraggio, E. Molho, A. Zaffaroni, On the notion of
proper efficiency in vector optimization, J. Optim. Theory Appl.
{\bf 82}(1994) 1-21.
\item{[7]} X. Y. Zheng, Proper efficiency in locally convex
 topological vector spaces, J. Optim. Theory  Appl. {\bf
 94} (1997)  469-486.
 \item{[8]} J. Liu, W. Song. On proper efficiencies in locally
 convex spaces -- a survey, Acta Math. Vietnam {\bf 26} (2001)
 301-312.
\item{[9]} J.  Horv\'{a}th, Topological Vector Spaces and
 Distributions, Vol.1, Addison-Wesley, Massachusetts,  1966.
\item{[10]}  K. D. Bierstedt, An introduction to locally convex
inductive limits, in: Functional Analysis and its Applications,
World Science, Singapore, 1988, pp. 35-133.
\item{[11]} J. H. Qiu, Weak property ($Y_0$) and regularity of inductive limits,
J. Math. Anal. Appl., {\bf 246} (2000)  379-389.
\item{[12]} J. H. Qiu, Y. Hao, Scalarization of Henig properly
efficient points in locally convex spaces, J. Optim. Theory Appl.
{\bf 147} (1010) 71-92.

\end{description}

\end{document}